\documentclass[11pt]{amsart} 
\usepackage{graphicx} 
\usepackage{amsmath, amssymb, amsthm}
\usepackage[total={6.5in,8.5in}]{geometry}
\usepackage[hidelinks]{hyperref}
\usepackage{xcolor}
\usepackage{todonotes}
\usepackage{comment}

\hypersetup{
    colorlinks=true,
    citecolor=blue,
    linkcolor=black,
    urlcolor=magenta
}

\definecolor{ObjectColor}{RGB}{0, 128, 0}
\definecolor{OptionColor}{RGB}{0, 128, 128}
\definecolor{FunctionColor}{RGB}{192, 0, 255}
\definecolor{TextColor}{RGB}{139, 34, 82}

\usepackage{listings}

\definecolor{ObjectColor}{RGB}{0, 128, 0}
\definecolor{OptionColor}{RGB}{0, 128, 128}
\definecolor{FunctionColor}{RGB}{192, 0, 255}
\definecolor{TextColor}{RGB}{139, 34, 82}

\lstdefinelanguage{Macaulay2}{
    sensitive=true,
    morecomment=[l]{--},
    morestring=[b]",
    morekeywords=[1]{ 
        loadPackage, 
        netList, 
        factor, 
        apply, 
        error,
        ideal,
        check,
        flatten,
        transpose,
        diff
    },
    morekeywords=[2]{ 
        elapsedTime,
        for,
        in,
        list,
        to,
        not,
        from
    },
    morekeywords=[3]{QQ, PolynomialRing}, 
    morekeywords=[4]{ 
        MonomialOrder,
        Lex,
        GRevLex,
        true,
        false,
        Depth,
        ReduceStrata,
        Verbose,
        CheckAssumption,
        Boxes,
    }
}

\theoremstyle{plain}
\newtheorem{theorem}{Theorem}[section]

\theoremstyle{definition}
\newtheorem{definition}[theorem]{Definition}
\newtheorem{example}[theorem]{Example}

\theoremstyle{remark}
\newtheorem{remark}[theorem]{Remark}

\newcommand{\QQ}{\mathbb Q}
\newcommand{\RR}{\mathbb R}

\title{Comprehensive Gr\"obner Bases in Macaulay2}

\author{María Agustina Cagliero}
\address{María Agustina Cagliero, KU Leuven, Belgium}
\email{mariaagustina.cagliero@kuleuven.be}

\author{Oliver Clarke}
\address{Oliver Clarke, Department of Mathematics, Durham University, UK}
\email{oliver.clarke@durham.ac.uk}

\author{Lorenzo De Biase}
\address{Lorenzo De Biase, ENEA Bologna Research Center, Italy}
\email{lorenzo.debiase@enea.it}

\author{Angelo El Saliby}
\address{Angelo El Saliby, Max Planck Institute for Mathematics in the Sciences, Leipzig, Germany}
\email{angelo.el.saliby@mis.mpg.de}

\author{Giulia Gaggero}
\address{Giulia Gaggero, McMaster University, Canada}
\email{gaggerog@mcmaster.ca}

\author{Weijia Wang}
\address{Weijia Wang, LIP6, Sorbonne Université, France}
\email{weijia.wang@lip6.fr}

\date{November 2026}

\begin{document}

\begin{abstract}
    Studying parameter-dependent polynomial systems often requires stratifying parameter space according to the algebraic or geometric properties of the specialised systems. Comprehensive Gröbner bases and systems, which provide Gröbner bases valid under every specialisation, either uniformly or on specified regions of parameter space, are classical tools for studying polynomial ideals under specialisation. We introduce \texttt{ComprehensiveGBs}, a \emph{Macaulay2} package for computing comprehensive Gröbner bases and systems, either over the full parameter space or under specified parameter conditions. We describe its use and implementation, with examples involving Weierstrass cubics and graphs.
\end{abstract}


\maketitle
\markboth{CAGLIERO, CLARKE, DE BIASE, EL SALIBY, GAGGERO, WANG}{COMPREHENSIVE GR\"OBNER BASES IN MACAULAY2}

\section{Introduction}
When the coefficients of a polynomial system depend on parameters, the Gr\"obner basis of the specialised system may change as the parameters vary,
and one would like to describe all of these Gröbner bases at once.
In 1992, Weispfenning introduced two notions for this purpose \cite{WEISPFENNING19921}.
A \emph{comprehensive Gr\"obner system} covers the parameter space by finitely many subsets and gives, for each subset, a set of polynomials that specialises to a Gr\"obner basis when the parameters are specialised to a point of that subset.
A \emph{comprehensive Gr\"obner basis} is a single finite set of polynomials whose specialisation is a Gr\"obner basis of the specialised ideal for every value of the parameters.
These notions have since been developed further by Montes~\cite{MONTES2002183},
by Suzuki and Sato~\cite{SS}, by Nabeshima~\cite{Nabeshima},
and by Kapur, Sun, and Wang~\cite{KAPUR201327}, among others.
For example, a single comprehensive Gröbner system determines which curves in the family of Weierstrass cubics $y^2=x^3+ax+b$ are singular, for all values of $a$ and $b$ simultaneously,
which we show in Example~\ref{example:weierstrass}.

Comprehensive Gröbner systems and bases are implemented in several computer algebra systems. The library \texttt{grobcov.lib} for Singular \cite{greuel2001singular,montes2012software} computes the Gröbner cover \cite{montes2010grobner}, a canonical comprehensive Gröbner system whose parameter sets are disjoint and locally closed. The package \texttt{PGB} for Risa/Asir \cite{noro1992risa,nabeshima2007pgb}, which grew out of Nabeshima's work \cite{Nabeshima}, computes comprehensive Gröbner systems and bases for polynomial rings, rings of differential operators, and modules. To our knowledge, no implementation was available in \emph{Macaulay2} \cite{M2}.

In this article, we introduce the package \texttt{ComprehensiveGBs} for \emph{Macaulay2}.
The package implements two complementary approaches to computing parametric Gröbner systems: the Suzuki--Sato algorithm \cite{SS} and the Kapur--Sun--Wang algorithm \cite{KAPUR201327}. The Suzuki--Sato algorithm computes a \emph{faithful} comprehensive Gröbner system, whose associated polynomial sets can be combined to form a comprehensive Gröbner basis. The Kapur--Sun--Wang algorithm typically produces fewer branches and requires less computation time than the Suzuki--Sato algorithm, although the resulting comprehensive Gröbner system need not be faithful. By providing both algorithms, the package allows users to choose between constructing a comprehensive Gröbner basis and efficiently computing Gröbner bases along with the corresponding conditions on the parameters.

The structure of the paper is as follows. In Section~\ref{section:background}, we fix our notation, recall definitions of comprehensive Gröbner systems and bases, and apply them to the family of Weierstrass cubics with our package. In Section~\ref{section:using-the-package}, we show how to use the package, including the options of each algorithm, with examples. Finally, in Section~\ref{section:design-choices}, we describe the data structures that the package uses internally, and other implementation choices.

\section{Background}\label{section:background}

Let $K$ be a field, and let $K[U]$ be the polynomial ring in the variables $ U=\{u_1,\ldots,u_m\}$, which we view as parameters. Let
$K[U][X]$, with $X\cap U=\emptyset$, be
the polynomial ring in the variables $X=\{x_1,\ldots,x_n\}$, with coefficients in the parameter ring.

Let $L$ be an algebraically closed extension field of $K$.
A \emph{specialisation} is a ring homomorphism
$\sigma:K[U]\to L$.
For $a=(a_1,\ldots,a_m)\in L^m$, we write
$\sigma_{a}$ for the specialisation given by
$\sigma_{a}(f)=f(a)$.
Every such specialisation extends naturally to a map
\[
\sigma_{\mathbf{a}}:K[U][X]\longrightarrow L[X]
\]
by applying $\sigma_{\mathbf{a}}$ coefficient-wise.

For $f_1,\ldots,f_k\in K[U]$, we denote by $V(f_1,\ldots,f_k)\subset L^m$ the affine variety of $f_1,\ldots,f_k$, and we extend the notation to subsets $F\subset K[U]$ likewise.

\textbf{Notation.} We fix monomial orders $\prec_X$ and $\prec_U$, and we equip $K[U,X]$ with the block order $\prec_{X,U}$, which compares the $X$-part first. All Gr\"obner bases in $L[X]$ are taken with respect to
$\prec_X$, unless otherwise stated. For a nonzero polynomial $g\in K[U][X]$, we write $\operatorname{lpp}_X(g)$ for its leading power product in $X$, and $\operatorname{lc}_X(g)\in K[U]$ for the corresponding coefficient.
Whenever we say that $\sigma_{a}(G)$ is a Gröbner basis, we understand that its zero elements are discarded, and we regard the empty set as a Gröbner basis of the zero ideal.

\begin{definition}
Let $E$ and $N$ be subsets of $K[U]$. The pair $(E,N)$ is called a \emph{parametric constraint}. It is called \emph{consistent} if $V(E)\setminus V(N)\neq\emptyset$, and \emph{inconsistent} otherwise.

Equivalently, $(E,N)$ is consistent if there exists at least one $f\in N$ such that
$f\notin\sqrt{\langle E\rangle}$.

 We say that the pair $(E,N)$ \emph{represents} the locally closed set $A=V(E)\setminus V(N)\subset L^m$.
\end{definition}

A parameter value belongs to $A=V(E)\setminus V(N)$ when every polynomial in $E$ vanishes on that value and the evaluation of at least one polynomial in $N$ is nonzero. For example, $N=\lbrace a,b\rbrace$ means $a\neq 0$ or $b\neq 0$, whereas $N=\lbrace ab\rbrace$ requires both $a\neq 0$ and $b\neq 0$. The choice $N=\lbrace 1\rbrace$ imposes no additional restriction.

\begin{definition}
Let $F$ be a subset of $K[U][X]$, let $A_1,\ldots,A_\ell$ be locally closed subsets of $L^m$,
let $G_1,\ldots,G_\ell$ be finite subsets of $K[U][X]$, and let $S$ be a subset of $L^m$ such that $S \subseteq A_1 \cup \cdots \cup A_\ell$. A finite set
\[
\mathcal{G} = \{(A_1,G_1),\ldots,(A_\ell,G_\ell)\}
\]
is called a \emph{comprehensive Gr\"obner system} on $S$ for $F$ if $\sigma_{\mathbf{a}}(G_i)$ is a Gr\"obner basis for the ideal $\langle \sigma_{\mathbf{a}}(F)\rangle$ in $L[X]$ for every $i=1,\ldots,\ell$ and every $\mathbf{a}\in A_i$.
Each pair $(A_i,G_i)$ is called a \emph{branch} of $\mathcal{G}$.

If $S=L^m$, then $\mathcal{G}$ is simply called a comprehensive Gr\"obner system for $F$.
\end{definition}

A comprehensive Gr\"obner system $\mathcal{G}=\{(A_1,G_1),\ldots,(A_\ell,G_\ell)\}$ for $F$ is called \emph{faithful} if, in addition, every element of $G_i$ belongs to $\langle F\rangle$ for $i=1,\ldots,\ell$.

\begin{definition}
Let $F$ be a subset of $K[U][X]$, and let $S$ be a subset of $L^m$.
A finite subset $G$ of $K[U][X]$ is called a \emph{comprehensive Gr\"obner basis} on $S$ for $F$ if $\sigma_{\mathbf{a}}(G)$ is a Gr\"obner basis for the ideal $\langle \sigma_{\mathbf{a}}(F)\rangle$ in $L[X]$ for every $\mathbf{a}\in S$.
If $S=L^m$, then $G$ is simply called a comprehensive Gr\"obner basis for $F$.
\end{definition}

A comprehensive Gr\"obner basis $G$ for $F$ is called \emph{faithful} if, in addition, every element of $G$ belongs to $\langle F\rangle$.

\begin{remark}
    Let $S$ and $F$ be as above. If $\mathcal{G}=\{(A_i,G_i)\}_{i=1}^{\ell}$ is a faithful comprehensive Gr\"obner system on $S$ for $F$, then $\bigcup_{i=1}^{\ell}G_i$ is a comprehensive Gr\"obner basis on $S$ for $F$.
\end{remark}

\begin{example}\label{example:weierstrass}
Consider the family of Weierstrass cubics
\[
C_{a,b}:y^2=x^3+ax+b,
\qquad (a,b)\in\mathbb{C}^2.
\]
We work over $K=\mathbb{Q}$, and we use a comprehensive
Gr\"obner system to study the singularities of the entire family simultaneously.

A point on $C_{a,b}$ is singular when both
partial derivatives vanish at it.
Thus, the singular points are the common zeros of the parametric ideal
\[
I=\left\langle
x^3+ax+b-y^2,\,-2y,\,3x^2+a
\right\rangle
\subseteq\mathbb{Q}[a,b][x,y].
\]
To compute a comprehensive Gr\"obner system for $I$ with
respect to the lexicographic order $x>y$, we define \texttt{F} as the list of the generators of $I$, and we call the function \texttt{comprehensiveGBMain} on it with the option \texttt{ReduceStrata} enabled.

\begin{lstlisting}
i1 : loadPackage "ComprehensiveGBs";
i2 : R = QQ[a,b][x,y, MonomialOrder => Lex];
i3 : F = {x^3 + a*x + b - y^2, -2*y, 3*x^2 + a};
i4 : G = comprehensiveGBMain(F, ReduceStrata => true);
i5 : netList G
     +------------------------------------------------------------------------+
     |         3      2      3      2                                         |
o5 = |({0}, {4a  + 27b }, {4a  + 27b })                                       |
     +------------------------------------------------------------------------+
     |    3      2             3      2              2               2        |
     |({4a  + 27b }, {a*b}, {4a  + 27b , y, 9b*x - 2a , 2a*x + 3b, 3x  + a})  |
     +------------------------------------------------------------------------+
     |                 3      2                     2                  2      |
     |({b, a}, {1}, {4a  + 27b , b*y, a*y, 9b*x - 2a , 2a*x + 3b, y, 3x  + a})|
     +------------------------------------------------------------------------+
\end{lstlisting}

The output is a comprehensive Gröbner system that consists of three branches. Each triple
$(E,N,G)$ represents the parameter region
$V(E)\setminus V(N)$, together with a set $G$ of polynomials whose
specialisation is a Gr\"obner basis on that region.

The first branch is given by $4a^3+27b^2\neq0$. After specialising $a$ and $b$ to such values, the polynomial $4a^3+27b^2$ becomes a nonzero constant and hence a unit. Thus, the specialised ideal is the unit ideal, whose reduced
Gr\"obner basis is $\{1\}$, so the singularity ideal has no solution. Hence, the curve $C_{a,b}$ is nonsingular.

The second branch is given by $4a^3+27b^2=0$ and $ab\neq0$. The specialised basis does not contain a unit, and therefore the curve has a singular point. From the polynomials $y$ and $2ax+3b$ in the displayed basis, we deduce that the point is $(-\frac{3b}{2a}, 0)$.

Finally, the third branch corresponds to $a=b=0$, for which the cubic $y^2=x^3$ has a singularity at the origin.

Notice that this comprehensive Gr\"obner system is
faithful: for each branch $(E_i,N_i,G_i)$, with
$i=1,2,3$, every polynomial in $G_i$ belongs to the
original ideal $\langle F\rangle$.
The first branch illustrates this.
Its basis is $G_1=\{4a^3+27b^2\}$, whose generator
belongs to $\langle F\rangle$.
Although $G_1^\prime=\{1\}$ would also be a valid
Gr\"obner basis on this region,
replacing $G_1$ by $G_1^\prime$ would destroy
faithfulness, since $1\notin\langle F\rangle$.
\end{example}

\section{Using the package}\label{section:using-the-package}

\subsection{The Suzuki--Sato algorithm.}

The function \texttt{comprehensiveGBMain} implements the Suzuki--Sato algorithm ~\cite{SS}, which computes a comprehensive Gr\"obner
system by recursively identifying parameter regions on which
a computed basis remains valid under specialisation.
We first demonstrate its use, and then we explain its available options.

\subsubsection{Usage.}\label{section:usage-example}
Let $F=\{ax+by\}$, where $a$ and $b$ are parameters and $x$ and $y$ are variables.
We distinguish parameters from variables with a tower of polynomial rings and equip it with a lexicographic order $x>y$.
\begin{lstlisting}
i1 : loadPackage "ComprehensiveGBs";
i2 : U = QQ[a,b];
i3 : R = U[x,y, MonomialOrder => Lex];
i4 : F = {a*x + b*y};
i5 : L = comprehensiveGBMain F;
i6 : netList L
     +------------------------------------+
o6 = |({0}, {a}, {a*x + b*y})             |
     +------------------------------------+
     |             2                      |
     |({a}, {b}, {a x + a*b*y, a*x + b*y})|
     +------------------------------------+
     |({b, a}, {1}, {a*x + b*y})          |
     +------------------------------------+
\end{lstlisting}

The function \texttt{comprehensiveGBMain} returns a list of triples.
Each triple $(E,N,G)$ describes the set $V(E)\setminus V(N)$, together with a list of polynomials $G$ whose specialisation to every point of that set is a Gröbner basis. In other words, a parameter value $(a,b) \in \overline \QQ^2$ lies in the set if for all $g \in E$ we have $g(a,b) = 0$ and there exists some $h \in N$ such that $h(a,b) \neq 0$.

With $f:= ax+by$, the three branches in this example are as follows.
\[
\begin{array}{c|c|c|c}
E & N & G & \text{conditions on the parameters} \\ \hline
\{0\}   & \{a\} & \{f\}    & a\ne 0 \\
\{a\}   & \{b\} & \{af,f\} & a=0\text{ and } b\ne 0 \\
\{a,b\} & \{1\} & \{f\}    & a=b=0.
\end{array}
\]
 The first branch corresponds to $a\ne 0$, where every specialisation of $f$ correspond to the unique element of a Gröbner basis of the specialised ideal; since $a$ is nonzero on this branch, it has leading monomial $x$.
 On the second branch, the polynomial $f$ specialises to the nonzero polynomial $by$, while $af$ specialises to zero. After discarding zero elements, the specialised Gröbner basis is therefore $\lbrace by\rbrace$.
On the last branch, $f$ specialise to zero. Then the specialised ideal is the zero ideal.

The function \texttt{comprehensiveGB} computes a comprehensive Gröbner basis instead, without any conditions on the parameters. The basis is obtained by combining the polynomial sets in the faithful Gröbner system computed by the function \texttt{comprehensiveGBMain}.
\begin{lstlisting}
i7 : comprehensiveGB F
                  2
o7 = {a*x + b*y, a x + a*b*y}
\end{lstlisting}
The result is a list of polynomials in $R$ not necessarily minimal.
For example, the set $\{ax + by\}$ alone is already a comprehensive Gr\"obner basis for $F$.

\subsubsection{Options.}
The option \texttt{Depth} bounds the recursion depth of
\texttt{comprehensiveGBMain}. Following from Section~\ref{section:usage-example}, the following counts the branches obtained with the depths $0$, $1$, and $2$.
\begin{lstlisting}
i8 : apply(0..2, d -> #comprehensiveGBMain(F, Depth => d))
o8 = (1, 2, 3)
\end{lstlisting}
With the depth $0$, the function returns only the branch with $a\ne 0$. The depth $1$ adds the branch with $a=0$ and $b\ne 0$ and the depth $2$ adds the branch at the origin.

By default, the value of the option \texttt{Depth} is $-1$, which imposes no depth bound.
With a nonnegative value, the function
may not explore the whole parameter space, so in general the result is a comprehensive Gröbner system on part of the parameter space only.
Likewise, with a nonnegative value of the option, the function \texttt{comprehensiveGB} is not guaranteed to return a comprehensive Gröbner basis on the full parameter space.

When the option \texttt{Verbose} is set to \texttt{true},
the function prints the lists $F$ and $S$
at each recursive call. The returned result does not change.
\begin{lstlisting}
i9 : comprehensiveGBMain(F, Verbose => true);
\end{lstlisting}

Recursive steps can lead to overlapping parameter regions.
The option \texttt{ReduceStrata} reduces this
redundancy by skipping components that are already
covered by the branches computed so far. The option does not, however, guarantee a disjoint or minimal Gröbner system. The default value of the option is false. In the following example, with two parameters and three variables, the option \texttt{ReduceStrata} reduces the number of returned branches from $230$ to $15$.

\begin{lstlisting}
i1 : loadPackage "ComprehensiveGBs";
i2 : U = QQ[a,b];
i3 : R = U[x,y,z, MonomialOrder => Lex];
i4 : F = {x^4 - a, y^5 - b, x + y - z};
i5 : L = comprehensiveGBMain(F, ReduceStrata => false);
i6 : Lred = comprehensiveGBMain(F, ReduceStrata => true);
i7 : (#L, #Lred)
o7 = (230,15)
\end{lstlisting}

The function \texttt{comprehensiveGBMain} also accepts a list $S$ of polynomials in $K[U]$ after the list $F$, in which case it computes a comprehensive Gröbner system on $V(S)$.
The list $S$ must satisfy
$V(S)\subseteq V\bigl(\langle F\rangle\cap K[U]\bigr)$,
outside of which the specialised ideal is the unit ideal.
By default, the function verifies this assumption before the computation starts, and it raises an error when the assumption fails.
The following code shows a case where the function raises an error.
\begin{lstlisting}[breaklines=true, breakatwhitespace=true]
i1 : loadPackage "ComprehensiveGBs";
i2 : U = QQ[a,b];
i3 : R = U[x, MonomialOrder => Lex];
i4 : F = {a*x, a*(a-1)};
i5 : Lrestricted = comprehensiveGBMain(F, {b});
error: V(S) is not contained in V(ideal F cap k[U]); pass CheckAssumption => false to skip this check
\end{lstlisting}

If the assumption is known to be true, then the option \texttt{CheckAssumption} can be set to \texttt{false} to skip this verification step. Otherwise, disabling the check may lead to an incorrect result. The following code restricts the computation of a previous example to $V(a)$.

\begin{lstlisting}
i1 : loadPackage "ComprehensiveGBs";
i2 : U = QQ[a,b];
i3 : R = U[x,y, MonomialOrder => Lex];
i4 : F = {a*x+b*y};
i5 : Lrestricted = comprehensiveGBMain(F, {a});
i6 : netList Lrestricted
     +------------------------------------+
     |             2                      |
o6 = |({a}, {b}, {a x + a*b*y, a*x + b*y})|
     +------------------------------------+
     |({b, a}, {1}, {a*x + b*y})          |
     +------------------------------------+
\end{lstlisting}
The computation now concerns only the cases with
$a=0$ and $b\neq 0$, and the case with $a=b=0$.
The branch with $a\ne 0$ lies outside the requested region.

\subsubsection{Comprehensive Gröbner Systems on Graphs.}
The function \texttt{comprehensiveGBOnGraph} constructs a
parametric polynomial system from a graph and
computes a comprehensive Gr\"obner system for it
using \texttt{comprehensiveGBMain} accepting all the same options.
Its input consists of a graph represented by
$\{V,E\}$ and a positive integer $d$, specifying
the dimension.

For each vertex $v$, the function introduces
coordinates $x_{v,1},\ldots,x_{v,d}$, which can be considered the position of that vertex, and for each
edge $\{i,j\}$, introduces a parameter $w_{ij}$ representing the squared Euclidean distance between its end-points.
The resulting polynomial system is
\[
    F=\left\{
    \sum_{k=1}^{d}(x_{i,k}-x_{j,k})^2-w_{ij}
    : \{i,j\}\in E
    \right\}
    \subseteq
    \mathbb{Q}[w_{ij}:\{i,j\}\in E]
    [x_{v,k}:v\in V,\ 1\leq k\leq d].
\]

The function returns a pair $(F, \mathbf G)$,
where $F$ is the polynomial system above and $\mathbf G$
is its comprehensive Gr\"obner system.
The graph can be supplied as a list of vertices
and a list of edge pairs, but also as an object
of type \texttt{Graph}. 
For example, consider the graph $G = (V, E)$ where
\[
V = \{1,2,3,4\} 
\quad \text{and} \quad
E = \{12, 13, 14, 24, 34\}.
\]
Let us set $d = 2$, i.e., we consider $G$ embedded in $\RR^2$ and let $F$ be the corresponding set of polynomials. We compute the generic stratum of a comprehensive Gr\"obner system for $F$ as follows.
\begin{lstlisting}
i1 : loadPackage "ComprehensiveGBs";
i2 : V = {1,2,3,4};
i3 : E = {{1,2},{1,3},{1,4},{2,4},{3,4}};
i4 : (F, GG) = comprehensiveGBOnGraph({V,E}, 2, ReduceStrata => true, Depth => 0);
i5 : netList F
     +--------------------------------------------------------+
     | 2      2                  2                  2         |
o5 = |x    + x    - 2x   x    + x    - 2x   x    + x    - w   |
     | 1,1    1,2     1,1 2,1    2,1     1,2 2,2    2,2    1,2|
     +--------------------------------------------------------+
     | 2      2                  2                  2         |
     |x    + x    - 2x   x    + x    - 2x   x    + x    - w   |
     | 1,1    1,2     1,1 3,1    3,1     1,2 3,2    3,2    1,3|
     +--------------------------------------------------------+
     | 2      2                  2                  2         |
     |x    + x    - 2x   x    + x    - 2x   x    + x    - w   |
     | 1,1    1,2     1,1 4,1    4,1     1,2 4,2    4,2    1,4|
     +--------------------------------------------------------+
     | 2      2                  2                  2         |
     |x    + x    - 2x   x    + x    - 2x   x    + x    - w   |
     | 2,1    2,2     2,1 4,1    4,1     2,2 4,2    4,2    2,4|
     +--------------------------------------------------------+
     | 2      2                  2                  2         |
     |x    + x    - 2x   x    + x    - 2x   x    + x    - w   |
     | 3,1    3,2     3,1 4,1    4,1     3,2 4,2    4,2    3,4|
     +--------------------------------------------------------+
i6 : netList for g in GG list {g_0, g_1}

     +---+-------------------------------------------------------------------+
     |   |                             2                  2                  |
o6 = |{0}|{w   w   w   w   w    - w   w   w   w    - w   w   w   w    + ... }|
     |   |  1,2 1,3 1,4 2,4 3,4    1,2 1,4 2,4 3,4    1,3 1,4 2,4 3,4        |
     +---+-------------------------------------------------------------------+
\end{lstlisting}
The generic stratum is the non-vanishing locus of the above polynomial, which can be factored as
\[
w_{3,4} w_{2,4} w_{1,4} (w_{1,3} - w_{1,4} - w_{3,4})(w_{1,2} - w_{1,4} - w_{2,4}).
\]
We may interpret this locus in terms of configurations of the graph, i.e., the parameters are in generic position when: none of the edges $34$, $24$, $14$ have length zero; the points $134$ are not colinear with $4$ lying between points $1$ and $3$; and the points $124$ are not colinear with $4$ lying between points $1$ and $2$.

\subsection{The Kapur--Sun--Wang algorithm.}

The function \texttt{parametricGBMain} implements the Kapur--Sun--Wang algorithm. This approach uses minimal Dickson bases to determine
the nonvanishing conditions needed for specialisation,
and it checks the consistency of the parameter conditions
during the computation~\cite{KAPUR201327}. We first introduce the definition of minimal Dickson basis of a system $H$, following the notation established in the background section. We then illustrate the use of the function and describe the available options.

\begin{definition}\cite[Definition~4.1]{KAPUR201327}
Fix a monomial order $\prec_X$, and let
$H \subset K[U][X]$ be a finite set of nonzero polynomials.
A subset $M \subseteq H$ is called a \emph{minimal Dickson basis}
of $H$ if the following two conditions hold.
\begin{enumerate}
    \item For every $h \in H$, there exists $g \in M$ such that
    $\operatorname{lpp}_X(g)$ divides $\operatorname{lpp}_X(h)$.
    \item For any two distinct elements $g$ and $g^\prime$ of $M$,
    $\operatorname{lpp}_X(g)$ does not divide
    $\operatorname{lpp}_X(g^\prime)$.
\end{enumerate}

A minimal Dickson basis need not be unique, since distinct
polynomials may have the same leading power product.
\end{definition}

\subsubsection{Usage.}
The function accepts a list $F$ of polynomials in $R=K[U][X]$, in which case it computes a comprehensive Gröbner system on the whole parameter space, as does the function \texttt{comprehensiveGBMain}.
The output is a list of lists $\{E_i,N_i,G_i\}$, where $G_i$ specialises to a Gr\"obner basis on $V(E_i)\setminus V(N_i)$.
The following example, with three parameters and two variables, illustrates the computation on the whole parameter space.
Since the branches are too wide to fit on one line, we display each of them as three rows labelled $E$, $N$, and $G$.
\begin{lstlisting}[
    label={ex:computation with lex}]
i1 : loadPackage "ComprehensiveGBs";
i2 : U = QQ[a,b,c];
i3 : R = U[x,y, MonomialOrder => Lex];
i4 : F = {a*x - b, b*y - a, c*x^2 - y, c*y^2 - x};
i5 : L = parametricGBMain F;
i6 : T = flatten apply(L, l -> transpose {{"E", "N", "G"}, l});
i7 : netList(T, Boxes => {{0, 3, 6, 9, 12}, {0, 1, 2}})
     +-+-------------------------------------------------+
o7 = |E|{0}                                              |
     | |    2         2       3     3   3     3   6    6 |
     |N|{b*c  - b, a*c  - a, b c - a , a c - b , a  - b }|
     |G|{1}                                              |
     +-+-------------------------------------------------+
     | |    2         2       3     3   3     3   6    6 |
     |E|{b*c  - b, a*c  - a, b c - a , a c - b , a  - b }|
     |N|{b}                                              |
     | |        2                                        |
     |G|{x - c*y , b*y - a}                              |
     +-+-------------------------------------------------+
     |E|{b, a}                                           |
     |N|{c}                                              |
     | |        2   3 4                                  |
     |G|{x - c*y , c y  - y}                             |
     +-+-------------------------------------------------+
     |E|{c, b, a}                                        |
     |N|{1}                                              |
     |G|{x, y}                                           |
     +-+-------------------------------------------------+
\end{lstlisting}

To restrict the computation, we specify the parameter
conditions together with the polynomial list.
The function \texttt{parametricGBMain} also accepts a list $E$ of polynomials in $K[U]$ and a list $N$ of polynomials in $K[U]$ before the list $F$, in which case it computes a comprehensive Gröbner system on $V(E)\setminus V(N)$.
In this interface, an empty list for $E$ is interpreted as $\lbrace 0\rbrace$ and an empty list for $N$ as $\lbrace 1\rbrace$, so empty input lists impose no parameter restrictions.
For instance, the following code restricts the computation for the same system to the parameter values with $a=b=0$ and
$c\ne 0$.
\begin{lstlisting}
i8 : Lrestricted = parametricGBMain({a, b}, {c}, F);
i9 : netList Lrestricted
     +------+---+--------------------+
     |      |   |        2   3 4     |
o9 = |{b, a}|{c}|{x - c*y , c y  - y}|
     +------+---+--------------------+
\end{lstlisting}
Here, the list $E=\{a,b\}$ gives the vanishing conditions, and the list $N=\{c\}$ gives the nonvanishing condition. The function \texttt{parametricGBMain} computes only on the requested region, and here, it returns the corresponding branch of the unrestricted computation.

\subsubsection{Minimal Dickson basis.}
The function \texttt{parametricGBMain} computes a minimal
Dickson basis internally to determine the
nonvanishing conditions used in the parameter
decomposition. This step is also available
independently through the function \texttt{minimalDicksonBasis}, which
takes a list of polynomials and selects a
minimal Dickson basis from it.
For example, consider the following list of polynomials.
\begin{lstlisting}
i1 : loadPackage "ComprehensiveGBs";
i2 : U = QQ[a, MonomialOrder => Lex];
i3 : R = U[x, y, MonomialOrder => Lex];
i4 : G = {a*x^2 - y, a*y^2 - 1, a*x - 1, (a + 1)*x - y, (a + 1)*y - a};
i5 : M = minimalDicksonBasis G;
\end{lstlisting}
The leading power products, taken with respect
to $x$ and $y$, are $x^2$, $y^2$, $x$, $x$, and $y$.
The minimal ones under divisibility are $x$ and
$y$, so the function retains one polynomial
for each of them. The choice need not be unique: both $ax-1$ and
$(a+1)x-y$ have leading power product $x$.
Among the candidates, the function first prefers leading coefficients shared by several candidates, then breaks ties using the degrees of their irreducible factors and other simplicity criteria, following the heuristics described
in~\cite[Section~7.1]{KAPUR201327}.
We explain this further in the design choices.
In this example, the function returns the following minimal Dickson basis.
\begin{lstlisting}
i6 : netList M
     +------------+
o6 = |(a + 1)x - y|
     +------------+
     |(a + 1)y - a|
     +------------+
\end{lstlisting}

\subsubsection{Options.} The option \texttt{Depth} bounds the recursion depth, with behavior identical to that of the function \texttt{comprehensiveGBMain}.
The option \texttt{ICheckLoops} affects the heuristic in a consistency check, which decides whether $f\in\sqrt{\langle E\rangle}$ for a constraint $(E,\lbrace f\rbrace)$.

If its value is a nonnegative integer $k$, then the check first squares $f$ up to $k$ times, and it falls back on the Rabinowitsch trick only if none of the powers $f^{2^i}$ is zero modulo $E$.
The default value \texttt{infinity} skips the squarings and applies the Rabinowitsch trick directly \cite[Chapter~4, Proposition~8]{cox2007ideals}.
The comparison below illustrates why the direct exact test can be preferable to repeated squaring.
We use the five-parameter, two-variable formulation of the Perspective-Three-Point (P3P) problem from computer vision given in \cite[Section~8]{KAPUR201327}.
The following code reports the measurements with \emph{Macaulay2} version 1.26.06 on a MacBook Pro (Apple M4 Pro, 2024) with 48 GB of RAM running macOS 27.2.1, which we use for all measurements presented in the paper.
\begin{lstlisting}
i1 : loadPackage "ComprehensiveGBs";
i2 : U = QQ[p, q, r, a, b, MonomialOrder => Lex];
i3 : R = U[x, y, MonomialOrder => GRevLex];
i4 : F = {(1 - a)*y^2 - a*x^2 - p*y + a*r*x*y + 1, (1 - b)*x^2 - b*y^2 - q*x + b*r*x*y + 1};
i5 : elapsedTime L = parametricGBMain F;
 -- .865066s elapsed
i6 : elapsedTime L4 = parametricGBMain(F, ICheckLoops => 4);
 -- 10.1557s elapsed
i7 : L === L4
o7 = true
\end{lstlisting}
With the value $k=5$, the same computation does not finish in $20$ minutes.

\subsection{Examples.}
We now apply both functions \texttt{comprehensiveGBMain} and \texttt{parametricGBMain} to several
parametric polynomial systems with various numbers
of parameters, variables, and generators. The monomial order is specified in each case.

\begin{example}[{\cite[Example~1]{SS}}]
We first consider the system $F=\{x^3-a,\ y^4-b,\ x+y-z\}$,
with two parameters and three variables.
We compute a comprehensive Gr\"obner system using
the function \texttt{comprehensiveGBMain} with the option \texttt{ReduceStrata} enabled.

\begin{lstlisting}  
i1 : loadPackage "ComprehensiveGBs";
i2 : R = QQ[a,b][x,y,z, MonomialOrder => Lex];
i3 : F = {x^3 - a, y^4 - b, x + y - z};
i4 : G = comprehensiveGBMain(F, ReduceStrata => true);
i5 : netList for g in G list {g_0, factor \ g_1}
     +------------------+----------------------------------------------------------+
     |                  |    4   2     4        3      4      3        4        3  |
o5 = |{0}               |{(b) (a) (729a  - 4096b )(729a  + 64b )(16767a  + 5632b )}|
     +------------------+----------------------------------------------------------+
     |                  |    3                                                     |
     |{b}               |{(a) }                                                    |
     +------------------+----------------------------------------------------------+
     |{b, a}            |{1}                                                       |
     +------------------+----------------------------------------------------------+
     |                  |    2                                                     |
     |{a}               |{(b) }                                                    |
     +------------------+----------------------------------------------------------+
     |     4        3   |    3   3                                                 |
     |{729a  - 4096b }  |{(b) (a) }                                                |
     +------------------+----------------------------------------------------------+
     |     4      3     |    3   2                                                 |
     |{729a  + 64b }    |{(b) (a) }                                                |
     +------------------+----------------------------------------------------------+
     |       4        3 |    4   3                                                 |
     |{16767a  + 5632b }|{(b) (a) }                                                |
     +------------------+----------------------------------------------------------+
\end{lstlisting}

The output consists of the parametric constraints of seven branches. The first branch applies outside the vanishing loci of the five irreducible factors, namely the coordinate axes and three curves. The remaining six branches cover their vanishing loci. The conditions specify where each polynomial set, which is not displayed in the output, remains a Gröbner basis. Notice that the polynomial sets in each branch depend on the chosen monomial order.
\end{example}

\begin{example}[{\cite[Example~3]{SS}}]
We next consider a system of five polynomials in
five variables, depending on four parameters.
The generators are constructed from
$f=ax_1^2+by_1$ and $g=cy_2^2+dx_2$ and their
partial derivatives, as shown below.
We apply both the function \texttt{comprehensiveGBMain}, with
the option \texttt{ReduceStrata} enabled, and the function \texttt{parametricGBMain},
and compare the elapsed times and the numbers of
returned branches.

\begin{lstlisting}
i1 : loadPackage "ComprehensiveGBs";
i2 : R = QQ[a,b,c,d][x_1,x_2,y_1,y_2,s, MonomialOrder => Lex];    
i3 : f = a*x_1^2 + b*y_1;       
i4 : g = c*y_2^2 + d*x_2;
i5 : F = {f, g, (x_1 - x_2)^2 + (y_1 - y_2)^2 - s, diff(x_1,f)*diff(y_2,g) - 
diff(y_1,f)*diff(x_2,g), diff(x_1,f)*(y_1-y_2) - diff(y_1,f)*(x_1 - x_2)};
i6 : elapsedTime GSS = comprehensiveGBMain(F, ReduceStrata => true);
 -- 3.17027s elapsed
i7 : elapsedTime GKSW = parametricGBMain F;
 -- .452324s elapsed
i8 : (#GSS, #GKSW)
o8 = (25, 17)
\end{lstlisting}

Here, we see that the Kapur--Sun--Wang algorithm produces fewer branches, and the computation requires much less time than that of the Suzuki--Sato algorithm.
\end{example}

\begin{example}
Now we show an example of \texttt{comprehensiveGBOnGraph}: consider the cycle on four vertices in dimension two.
We display the total number of branches and the
parameter conditions of the first three.
\begin{lstlisting}
i1 : loadPackage "ComprehensiveGBs";
i2 : V = {1,2,3,4};
i3 : E = {(1,2),(1,3),(3,4),(2,4)};
i4 : G = {V,E};
i5 : (F, GG) = comprehensiveGBOnGraph(G, 2);
i6 : #GG
o6 = 11
i7 : netList for j from 0 to 2 list (GG_j)_{0,1}
     +---------------------------------------+
o7 = |({0}, {w   })                          |
     |        3,4                            |
     +---------------------------------------+
     |           2                  2        |
     |({w   }, {w   w   w    - w   w   w   })|
     |   3,4     1,2 1,3 2,4    1,2 1,3 2,4  |
     +---------------------------------------+
     |                 2              2      |
     |({w   , w   }, {w   w    - w   w   })  |
     |   2,4   3,4     1,2 1,3    1,2 1,3    |
     +---------------------------------------+
\end{lstlisting}

The first branch applies whenever $w_{3,4}\neq 0$. The last two branches are on $w_{3,4}=0$ and distinguish further vanishing conditions on the parameters. Hence, the output provides a comprehensive Gröbner system on the whole parameter space.
\end{example}

\begin{example}
We finally return to the example in Section~\ref{ex:computation with lex}, now using graded
reverse lexicographic order on $x,y$.
\begin{lstlisting}
i1 : loadPackage "ComprehensiveGBs";
i2 : U = QQ[a,b,c, MonomialOrder => GRevLex];
i3 : R = U[x,y, MonomialOrder => GRevLex];
i4 : F = {a*x - b, b*y - a, c*x^2 - y, c*y^2 - x};
i5 : L = parametricGBMain F;
i6 : T = flatten apply(L, l -> transpose {{"E", "N", "G"}, l});
i7 : netList(T, Boxes => {{0, 3, 6, 9, 12}, {0, 1, 2}})
     +-+-------------------------------------------------+
o7 = |E|{0}                                              |
     | |    2         2       3     3   3     3   6    6 |
     |N|{b*c  - b, a*c  - a, b c - a , a c - b , a  - b }|
     |G|{1}                                              |
     +-+-------------------------------------------------+
     | |    2         2       3     3   3     3   6    6 |
     |E|{b*c  - b, a*c  - a, b c - a , a c - b , a  - b }|
     |N|{b}                                              |
     |G|{b*x - a*c*y, b*y - a}                           |
     +-+-------------------------------------------------+
     |E|{b, a}                                           |
     |N|{c}                                              |
     | |    2         2                                  |
     |G|{c*x  - y, c*y  - x}                             |
     +-+-------------------------------------------------+
     |E|{c, b, a}                                        |
     |N|{1}                                              |
     |G|{x, y}                                           |
     +-+-------------------------------------------------+
\end{lstlisting}

The parameter conditions remain unchanged, while the bases of the second and third branches differ. The first and fourth bases remain $\lbrace 1\rbrace$ and $\lbrace x,y\rbrace$. This example also illustrates that the comprehensive
Gr\"obner system returned by \texttt{parametricGBMain}
need not be faithful.
Indeed, the first branch has polynomial basis $\{1\}$,
but $1$ does not belong to the original ideal.
\end{example}

\section{Design choices}\label{section:design-choices}

\subsection{Rings and maps}

Both algorithms pass polynomials between several rings, whose monomial orders are compatible with that of $R=K[U][X]$.
The flattened ring $K[X,U]$ makes the parameters available as variables for Gröbner basis and elimination computations.
The Suzuki--Sato algorithm also uses an auxiliary variable $\ell$ for preserving faithfulness during recursion.

\begin{table}[htbp]
    \centering
    \begin{minipage}[t]{0.46\textwidth}
        \vspace{0pt}
        \centering
        \small
        \setlength{\tabcolsep}{4pt}
        \begin{tabular}{lll}
            \hline
            Name & Ring & Order \\
            \hline
            \texttt{KU}
                & $K[U]$ & $U$ \\
            \texttt{RExt}
                & $K[\ell,X,U]$ & $\ell;\ X;\ U$ \\
            \texttt{RFlat}
                & $K[X,U]$ & $X;\ U$ \\
            \texttt{RExt'}
                & $K[U][\ell,X]$ & $\ell;\ X$ \\
            \texttt{RFlatl}
                & $K[X,U][\ell]$ & $\ell$ \\
            \hline
        \end{tabular}
        \caption{Rings constructed by
        \texttt{cgbDataFromRings}.}
        \label{tab:cgb-rings}
    \end{minipage}\hfill%
    \begin{minipage}[t]{0.52\textwidth}
        \vspace{0pt}
        \centering
        \small
        \setlength{\tabcolsep}{4pt}
        \begin{tabular}{ll}
            \hline
            Name & Map \\
            \hline
            \texttt{RtoRExt}
                & $R \to \texttt{RExt}$ \\
            \texttt{RExttoRFlatl}
                & $\texttt{RExt} \to \texttt{RFlatl}$ \\
            \texttt{RExttoRExt'}
                & $\texttt{RExt} \to \texttt{RExt'}$ \\
            \texttt{RExttoR}
                & $\texttt{RExt} \to R$ \\
            \texttt{KUtoRFlat}
                & $\texttt{KU} \to \texttt{RFlat}$ \\
            \texttt{RFlattoR}
                & $\texttt{RFlat} \to R$ \\
            \texttt{KUtoR}
                & $\texttt{KU} \to R$ \\
            \texttt{RtoRFlat}
                & $R \to \texttt{RFlat}$ \\
            \texttt{RFlattoKU}
                & $\texttt{RFlat} \to \texttt{KU}$ \\
            \hline
        \end{tabular}
        \caption{Ring maps constructed by
        \texttt{cgbDataFromRings}.}
        \label{tab:cgb-maps}
    \end{minipage}
\end{table}

The internal function \texttt{cgbDataFromRings} constructs these rings and the ring maps between them, which are listed in Table~\ref{tab:cgb-rings} and Table~\ref{tab:cgb-maps}.
The function creates the variables of the rings from local symbols, so it never overwrites a symbol of the user. Moreover, the function stores the rings and maps in a hash table, which is cached in $R$.

\subsection{The Suzuki--Sato algorithm}

Our implementation follows Algorithm
\texttt{CGBMain} of \cite{SS}, with the following two changes.
At each recursive step, let $S$ be the input parameter equations. Let $h_1,\ldots,h_t\in K[U]$ be the leading coefficients selected by Algorithm \texttt{CGBMain}, and denote $h=\operatorname{lcm}(h_1,\ldots,h_t)$.

\begin{itemize}
    \item The algorithm recurses on $S\cup\lbrace h^\prime_i\rbrace$ for the irreducible factors $h^\prime_i$ of $h$, rather than on $S\cup\lbrace h_i\rbrace$. Both choices cover $V(h)=\bigcup_i V(h_i)=\bigcup_i V(h^\prime_i)$, while factorisation replaces the coefficient equations by their irreducible factors before recursion.
    \item When the option \texttt{ReduceStrata} is set to \texttt{true}, the algorithm does not recurse on $S\cup\lbrace h^\prime_i\rbrace$ directly.
    Instead, it computes the minimal primes $P$ of $\langle S\cup\lbrace h^\prime_i\rbrace\rangle$ with the package \texttt{MinimalPrimes} \cite{MinimalPrimesSource},
    and it recurses on each $P$ only if $V(P)$ is not covered by the branches $A_j=V(E_j)\setminus V(N_j)$ computed so far. \\
    The coverage test is analogous to that of \cite[Table 6]{montes2010grobner}. It takes an ideal $I$ as input and returns true if $I$ is the unit ideal. Otherwise, for each minimal prime $Q$ of $I$, it looks for an index $j$ such that $E_j\subset P$ and $N_j\not\subset P$. If none exists, then the test returns false. Otherwise, $V(P)\setminus V(N_i)$ is covered by $A_i$, and the remaining set is tested recursively with input $P+\langle N_i\rangle$, with the test returning true if and only if all the recursive tests succeed.
\end{itemize}

\subsection{The Kapur--Sun--Wang algorithm}

Our implementation follows \cite[Section 4.1]{KAPUR201327} with two changes, both of which the paper mentioned outside its pseudocode.
At each recursive step, let $G$ be the computed Gröbner basis in $K[X,U]$. Let $G_r$ be its parameter-only part, and let $G_m$ be a minimal Dickson basis of the remaining polynomials, viewed in $K[U][X]$.
\begin{itemize}
    \item The branching uses the irreducible factors $h^\prime_i$ of the leading coefficients of the polynomials in $G_m$, as for the implementation of the Suzuki--Sato algorithm.
    \item The recursive calls receive Gröbner bases of the ideals $\langle G_r\cup\lbrace h^\prime_i\rbrace\rangle$ instead of the lists $G_r\cup\lbrace h^\prime_i\rbrace$. This is also what the authors do in \cite[Example 6.1]{KAPUR201327}, and the consistency checks rely on such bases.
\end{itemize}

\subsubsection*{Choosing a minimal Dickson basis}

The function \texttt{minimalDicksonBasis} sorts the candidates with a minimal leading power product, and then it keeps the first polynomial for each leading product. The sorting uses heuristics from \cite[Section 7.1]{KAPUR201327}. It prefers a polynomial whose leading coefficient occurs most often among the candidates with a minimal leading power product. If there are ties, then it prefers a polynomial with a ``simpler'' leading coefficient, i.e., one whose irreducible factors have the smallest maximum degree, then one whose leading coefficient has fewer terms, and then one whose leading coefficient has lower degree. Any remaining ties are broken by comparing the string representations of the polynomials, which makes the result deterministic. This tie-breaking order can affect the factors used for branching and hence the scale of the computation.

\subsubsection*{Consistency checks}

We check the consistency of a constraint $(E,\lbrace f\rbrace)$, where $E$ is a Gröbner basis of $\langle E\rangle$, with the checks of \cite[Section~5.3]{KAPUR201327}. The checks are with the same name and order as for \cite{KAPUR201327}: the normal form of $f$ modulo $E$, then the internal function \texttt{zeroDimCheck} if $\dim(\langle E\rangle)=0$, or the internal heuristic functions \texttt{cCheck} and \texttt{iCheck} otherwise.
Nevertheless, we depart from the paper in the last check. With the default value $k=\infty$ of the option \texttt{ICheckLoops}, the function \texttt{iCheck} skips the heuristics and tests whether $f\in\sqrt{\langle E\rangle}$ with the Rabinowitsch trick, applied to the square-free part of $f$.

For a constraint $(E,N)$ with several polynomials in $N$, we run these checks on each $f\in N$, and deduce that $(E,N)$ is consistent as soon as one of the constraints $(E,\lbrace f\rbrace)$ is consistent. If $k$ is finite and some of these checks are inconclusive, then we fall back on the Rabinowitsch trick for the polynomials in $N$. During consistency checks, heuristics within the internal functions could affect the computation path, without, however, affecting the final result.

\section*{Acknowledgements}

We would like to thank Diane Maclagan for organising the 2026 \textit{Warwick Macaulay2 Workshop} where this project started.
We would also like to thank Woody Cohen and Sam Knight for their valuable contributions to the development of the \texttt{ComprehensiveGBs} package.

María Agustina Cagliero is partially supported by the FWO grants G0F5921N (Odysseus) and G023721N, and by the KU Leuven grant iBOF/23/064.
Lorenzo De Biase is partially supported by the Programme Agreement between the Italian Ministry of the Environment and Energy Security  and ENEA for Electric System Research, Project 1.7, WP1, LA2.8 (2025–2027).
Giulia Gaggero is funded by SNF Postdoc.Mobility grant P500PT\textunderscore225436.
Weijia Wang is supported by the ANR Project ANR-22-CE91-0007 ``EAGLES''
and by grant FA8655-25-1-7469 of the European Office of Aerospace Research and Development of the Air Force Office of Scientific Research (AFOSR).

\subsection*{Use of AI} Large language models (LLM) from \emph{OpenAI} and \emph{Anthropic} were used only for lookups, explanations, or advice. They served as enhanced search engines to help us understand the referenced literature. No LLM was used to write any code for our program or any text in this paper.

\section*{Supplement}

The \href{https://github.com/Macaulay2/Workshop-2026-Warwick/tree/CGBs}{online supplement} includes the Git repository for the package \texttt{ComprehensiveGBs}. This paper is based on version 1.0 of the package, at commit \texttt{f6b6ae6}.

\bibliographystyle{amsalpha}
\bibliography{references}

@Misc{M2,
author = {Grayson, Daniel R. and Stillman, Michael E.},
title = {Macaulay2, a software system for research in algebraic geometry},
howpublished = {Available at \url{http://www2.macaulay2.com}}
}

@misc{MinimalPrimesSource,
title = {{MinimalPrimes: minimal primes and radical routines for ideals. Version~0.10}},
author = {Frank Moore and Mike Stillman and Franziska Hinkelmann and Justin Chen and Mahrud Sayrafi},
howpublished = {A \emph{Macaulay2} package available at
  \url{https://github.com/Macaulay2/M2/tree/stable/M2/Macaulay2/packages}}
}

@article{montes2010grobner,
  title={Gr{\"o}bner bases for polynomial systems with parameters},
  author={Montes, Antonio and Wibmer, Michael},
  journal={Journal of Symbolic Computation},
  volume={45},
  number={12},
  pages={1391--1425},
  year={2010},
  publisher={Elsevier}
}

@article{montes2012software,
  title={Software for computing the Gr{\"o}bner cover of a parametric ideal},
  author={Montes, Antonio},
  journal={ACM Communications in Computer Algebra},
  volume={45},
  number={3/4},
  pages={180--182},
  year={2012},
  publisher={ACM, New York, NY, USA}
}

@inproceedings{noro1992risa,
  title={Risa/Asir—a computer algebra system},
  author={Noro, Masayuki and Takeshima, Taku},
  booktitle={Papers from the international symposium on Symbolic and algebraic computation},
  pages={387--396},
  year={1992}
}

@article{greuel2001singular,
  title={Singular, a computer algebra system for polynomial computations},
  author={Greuel, Gert-Martin and Pfister, Gerhard and Sch{\"o}nemann, Hans},
  journal={Symbolic computation and automated reasoning},
  pages={227--233},
  year={2001}
}

@article{WEISPFENNING19921,
title = {Comprehensive Gröbner bases},
journal = {Journal of Symbolic Computation},
volume = {14},
number = {1},
pages = {1-29},
year = {1992},
issn = {0747-7171},
doi = {https://doi.org/10.1016/0747-7171(92)90023-W},
url = {https://www.sciencedirect.com/science/article/pii/074771719290023W},
author = {Volker Weispfenning}
}

@article{KAPUR201327,
title = {An efficient algorithm for computing a comprehensive Gröbner system of a parametric polynomial system},
journal = {Journal of Symbolic Computation},
volume = {49},
pages = {27-44},
year = {2013},
note = {The International Symposium on Symbolic and Algebraic Computation},
issn = {0747-7171},
doi = {https://doi.org/10.1016/j.jsc.2011.12.015},
url = {https://www.sciencedirect.com/science/article/pii/S0747717111002082},
author = {Deepak Kapur and Yao Sun and Dingkang Wang}
}

@inproceedings{SS,
author = {Suzuki, Akira and Sato, Yosuke},
title = {A simple algorithm to compute comprehensive Gr{\"o}bner bases using Gr{\"o}bner bases},
year = {2006},
isbn = {1595932763},
publisher = {Association for Computing Machinery},
address = {New York, NY, USA},
url = {https://doi.org/10.1145/1145768.1145821},
doi = {10.1145/1145768.1145821},
booktitle = {Proceedings of the 2006 International Symposium on Symbolic and Algebraic Computation},
pages = {326–331},
numpages = {6},
location = {Genoa, Italy},
series = {ISSAC '06}
}

@article{MONTES2002183,
title = {A New Algorithm for Discussing Gröbner Bases with Parameters},
journal = {Journal of Symbolic Computation},
volume = {33},
number = {2},
pages = {183-208},
year = {2002},
issn = {0747-7171},
doi = {https://doi.org/10.1006/jsco.2001.0504},
url = {https://www.sciencedirect.com/science/article/pii/S0747717101905048},
author = {Antonio Montes}
}

@article{Nabeshima,
author = {Nabeshima, Katsusuke},
year = {2007},
month = {07},
pages = {299-306},
title = {A speed-up of the algorithm for computing comprehensive Gröbner systems},
journal = {Proceedings of the International Symposium on Symbolic and Algebraic Computation, ISSAC},
doi = {10.1145/1277548.1277589}
}

@article{nabeshima2007pgb,
author = {Nabeshima, Katsusuke},
title = {PGB: a package for computing parametric Gr{\"o}bner and related objects},
year = {2007},
issue_date = {September 2007 issue 161},
publisher = {Association for Computing Machinery},
address = {New York, NY, USA},
volume = {41},
number = {3},
issn = {1932-2232},
url = {https://doi.org/10.1145/1358190.1358198},
doi = {10.1145/1358190.1358198},
journal = {ACM Commun. Comput. Algebra},
month = sep,
pages = {104–105},
numpages = {2}
}

@book{cox2007ideals,
  title={Ideals, varieties, and algorithms: an introduction to computational algebraic geometry and commutative algebra},
  author={Cox, David A and Little, John B and O'shea, Donal and Little, John},
  volume={10},
  year={2007},
  publisher={Springer}
}

\end{document}